\documentclass[11pt]{article}

\usepackage[cmtip,arrow]{xy}
\usepackage{amsmath,amssymb,enumerate,pb-diagram,pb-xy}
\def\today{\number\day .\number\month .\number\year}
\usepackage[applemac]{inputenc}

\def\1{{\bf 1}}
\def\A{{\mathbb A}}

\def\al{\alpha}
\def\bs{\backslash}
\def\C{{\mathbb C}}

\def\CF{{\cal F}}
\def\CG{{\cal G}}
\def\CJ{{\cal J}}

\def\CV{{\cal V}}
\def\comm{\operatorname{comm}}
\def\df{\ \stackrel{\mbox{\rm\tiny def}}{=}\ }
\def\e{\emph}

\def\Ext{\operatorname{Ext}}
\def\fin{{\rm fin}}
\def\g{{\mathfrak g}}
\def\G{{\mathbb G}}
\def\ga{\gamma}
\def\Ga{\Gamma}

\def\H{\operatorname H}

\def\Hom{{\rm Hom}}

\def\k{{\mathfrak k}}

\def\La{\Lambda}

\def\loc{{\rm loc}}
\def\mathqed{\vspace{-30pt}\\ \qed}

\def\Mod{\operatorname{Mod}}
\def\N{{\mathbb N}}

\def\p{{\mathfrak p}}

\def\ph{\varphi}
\def\prf{{\bf Proof: }}

\def\Q{{\mathbb Q}}
\def\qed{\ifmmode\eqno \square 
		\else\noproof\vskip 12pt plus 3pt minus 9pt \fi}
\def\noproof{{\unskip\nobreak\hfill\penalty50\hskip2em\hbox{}%
     \nobreak\hfill $\square$\parfillskip=0pt%
     \finalhyphendemerits=0\par}}
\def\R{{\mathbb R}}

\def\res{\operatorname{res}}
\def\Sheaf{\operatorname{Sheaf}}

\def\sm{\smallsetminus}

\def\what{\widehat}
\def\Z{{\mathbb Z}}

\def\({\left(}
\def\){\right)}
\def\={{\ =\ }}

\newcommand{\tto}[1]{\stackrel{#1}{\longrightarrow}}

\newcommand{\ol}[1]{\overline{#1}}
\newcommand{\stack}[2]{\genfrac{}{}{0pt}{1}{#1}{#2}}

\newtheorem{theorem}{Theorem}[section]

\newtheorem{lemma}[theorem]{Lemma}

\newtheorem{proposition}[theorem]{Proposition}

\newtheorem{defi}[theorem]{Definition}
\newenvironment{definition}[0]{\begin{defi}\rm}
{\end{defi}}
\newtheorem{exmples}[theorem]{Examples}

\begin{document}

\pagestyle{myheadings} \markright{HIGHER ORDER INVARIANTS}

\title{Higher order invariants in the case of compact quotients}
\author{Anton Deitmar\\ \ \\ 
Central European Journal of Mathematics 9.1, 85-101 (2011)}
\date{}
\maketitle

{\bf Abstract:}
We present the theory of higher order invariants and higher order automorphic forms in the simplest case, that of a compact quotient.
In this case many things simplify and we are thus able to prove a more precise structure theorem than in the general case.

\tableofcontents

\section*{Introduction}
Higher order modular forms show up in various contexts, for instance in percolation theory \cite{KZ}, or in the theory of Eisenstein-series twisted by modular symbols \cite{CDO,GG,Gold}.
Finally, spaces of higher order forms are natural receptacles of converse theorems \cite{F,IM}.

L-functions of second order forms have been studied in \cite{DKMO}, Poincar\'e series attached to higher order forms have been investigated in \cite{IO}, dimensions of spaces of second order forms have been determined in \cite{DO,DiamSim}.
Higher order cohomology has been introduced and an Eichler-Shimura type theorem has been proven in \cite{ES}.
In \cite{HOAut} a program has been started which aims at an understanding of the theory of higher order forms from an representation-theoretical point of view.

The paper \cite{HOinv} contains a structure theorem showing that spaces of higher order automorphic forms are in a natural way subspaces of tensor products of automorphic representation spaces.
In this paper we restrict to the simple case of a compact quotient in which case we are able to 
\begin{itemize}
\item remove the ``subspace'' part from the assertion, i.e., get a more precise statement on the structure of higher order forms, and
\item extend the theory to $L^2$-invariants instead of smooth invariants only. 
\end{itemize}
Also, this presentation has the advantage of being much simpler than in the general setting, for instance, in the case of compact quotient, there is no need to factor out cuspidal ideals and so all formulae look simpler.

\section{Higher order invariants}
Let $R$ be a commutative ring with unit.
In the main applications, $R$ will be the field of complex numbers.
Let $\Ga$ be a group and let $I$ denote the augmentation ideal in the group algebra $A=R[\Ga]$.
For an $A$-module $V$ we define the $R$-module
$$
\H_{q}^0(\Ga,V)
$$
of invariants of order $q\ge 1$ to be the set of all $v\in V$ with $I^{q}v=0$.
Note that for $q=1$ one gets the usual invariants
$$
\H_{1}^0(\Ga,V)\= \H^0(\Ga,V)\= V^\Ga.
$$
There is a natural identification
$$
\H_{q}^0(\Ga,V)\ \cong\ \Hom_A(A/I^{q},V).
$$

The sets $\H_{q-1}^0(\Ga,V)\subset \H_{q}^0(\Ga,V)$ form a filtration on $V$ which is not necessarily exhaustive.
Let
$$
\ol \H_{q}^0(\Ga,V)\df \H_{q}^0(\Ga,V)/\H_{q-1}^0(\Ga,V)
$$
be the $q$-th graded piece, where we allow $q=1,2,\dots$ by formally setting $\H_{0}^0(\Ga,V)=0$.

For an abelian group $G$, let $\Hom(\Ga,G)$ denote the set of all group homomorphisms $\Ga\to G$.

\begin{lemma}\label{1.1*}
Suppose that $R$ is a field.
If $\Ga$ is finitely generated as a group or if $W$ is finite-dimensional, then
$$
\Hom(\Ga,W)\cong \Hom(\Ga,R)\otimes W,
$$
where the tensor product is over $R$.
\end{lemma}

\prf
This is easy to see.
\qed

\begin{definition}
We introduce the \e{order-lowering-homomorphism}
$$
\La\=\La_{q}: \ol \H_{q}^0(\Ga,V)\to\Hom(\Ga,\ol \H_{q-1}^0(\Ga,V))
$$
given by
$$
\La(v)(\ga)\=(\ga-1)v.
$$
To see that $\psi(v)$ is indeed a group homomorphism, note that in the group algebra $\C[\Ga]$ one has $(\tau\ga-1)\equiv (\tau-1)+(\ga-1)\mod I^2$.
\end{definition}

\begin{proposition}
\begin{enumerate}[\rm (a)]
\item The order-lowering operator $\La$ is injective.
\item
If $R$ is a field and $\Hom(\Ga,R)=\H^1(\Ga,R)=0$, then there are no non-trivial higher order invariants, i.e., one has
$$
\H_{q}^0(\Ga,V)\= V^\Ga
$$
for all $q\ge 1$.
\item
If $R$ is a field and $\Ga$ is finitely generated or $\dim V<\infty$, then the order-lowering operator induces an injective $\Ga$-equivariant linear map
$$
\La^q: \ol \H_{q}^0(\Ga,V)\hookrightarrow \Hom(\Ga,R)^{\otimes (q-1)}\otimes V^\Ga.
$$

\end{enumerate}\end{proposition}

\prf
Part (a) is clear.
Part (b) becomes clear by choosing a basis of $\ol H_{q-1}^0(\Ga,V)$.
The last assertion follows from Lemma \ref{1.1*}.
\qed

\begin{lemma}
Let $V$ be an $R[\Ga]$-module which is torsion-free as $\Z$-module.
Let $\Sigma\subset\Ga$ be a subgroup of finite index.
Then the natural restriction map
$$
\res: \ol \H_{q}^0(\Ga,V)^0\ \to\ \ol \H_{q}^0(\Sigma,V)
$$
is injective.
\end{lemma}

Note that the torsion-condition is automatic if $R$ contains the field $\Q$.

\prf
By induction on $q$.
\qed

\section{Higher order cohomology}
For an $R[\Ga]$-module $V$ we define the \e{higher order cohomology} to be
$$
\H_{q}^p(\Ga,V)\=\Ext^p_A(A/I^{q},V).
$$
For $q=1$ this is ordinary group cohomology.

For an $R$-module $M$ and a set $S$ we write $M^S$ for the $R$-module of all maps from $S$ to $M$.
Then $M^\emptyset$ is the trivial module $0$. 
Up to isomorphy, the module $M^S$ depends only on the cardinality of $S$.
It therefore makes sense to define $M^c$ for any cardinal number $c$ in this way.
Note that $I^q/I^{q+1}$ is a free $R$-module.
Define
$$
N_{\Ga}(q)\df \dim_R I^q/I^{q+1}.
$$
Then $N_{\Ga}(q)$ is a possibly infinite cardinal number.

\begin{lemma}\label{1.1}
\begin{enumerate}[\rm (a)]
\item For every $q\ge 1$ there is a natural exact sequence
\begin{multline*}
0\to H_{q}^0(\Ga,V)\to H_{q+1}^0(\Ga,V)\to H^0(\Ga,V)^{N_{\Ga}(q)}\to\\
\to H_{q}^1(\Ga,V)\to H_{q+1}^1(\Ga,V)\to H^1(\Ga,V)^{N_{\Ga}(q)}\to\dots\\
\dots\to H_{q}^p(\Ga,V)\to H_{q+1}^p(\Ga,V)\to H^p(\Ga,V)^{N_{\Ga}(q)}\to\dots
\end{multline*}
\item Suppose that for a given $p\ge 0$ one has $H^p(\Ga,V)=0$.
Then $H_{q}^p(\Ga,V)=0$ for every $q\ge 1$.
In particular, if $V$ is acyclic as $\Ga$-module, then $H_{q}^p(\Ga,V)=0$ for all $p\ge 1$, $q\ge 1$.
\end{enumerate}\end{lemma}

Note that $\H^p(\Ga,V)^{N_{\Ga}(q)}$ here comes up as
$\Ext_A^p(I^q/I^{q+1},V)$,
which is important for the functoriality in the group $\Ga$  of the above sequence.

\prf
Consider the exact sequence
$$
0\to I^q/I^{q+1}\to A/I^{q+1}\to A/I^q\to 0.
$$
As an $A$-module, $I^q/I^{q+1}$ is isomorphic to a direct sum $\bigoplus_\al R_\al$ of copies of $R=A/I$.
So we conclude that for every $p\ge 0$,
$$
\Ext_A^p(I^q/I^{q+1},V)\ \cong\ \prod_\al\Ext_A^p(R,V)\ \cong\ H^p(\Ga,V)^{N_{\Ga,\Sigma}(q)}.
$$
The long exact $\Ext$-sequence induced by the above short sequence is
\begin{multline*}
0\to\Hom_A(A/I^q,V)\to\Hom_A(A/I^{q+1},V)\to \Hom_A(I^q/I^{q+1},V)\to\\
\to\Ext_A^1(A/I^q,V)\to\Ext_A^1(A/I^{q+1},V)\to \Ext_A^1(I^q/I^{q+1},V)\to\\
\to\Ext_A^2(A/I^q,V)\to\Ext_A^2(A/I^{q+1},V)\to \Ext_A^2(I^q/I^{q+1},V)\to\dots
\end{multline*}
This is the claim (a).
For (b) we proceed by induction on $q$.
For $q=1$ the claim follows from $H^p_1(\Ga,V)=H^p(\Ga,V)$.
Inductively, assume the claim proven for $q-1$ and $H^p(\Ga,V)=0$.
As part of the above exact sequence, we have the exactness of
$$
H_{q-1}^p(\Ga,V)\to H_{q}^p(\Ga,V)\to H^p(\Ga,V)^{N_{\Ga}(q)}.
$$
By assumption, we have $H^p(\Ga,V)^{N_{\Ga,\Sigma}(q)}=0$ and by induction hypothesis the module $H_{q-1}^p(\Ga,V)$ vanishes.
This implies $H_{q}^p(\Ga,V)=0$ as well.
\qed

\begin{lemma}\label{Prop2.1.3}
For every $q\ge 1$ there is a natural isomorphism
$$
\Hom_R(I^q/I^{q+1},R)\ \cong\ \H_{q}^1(\Ga,R).
$$
\end{lemma}

\prf
The exact sequence
$$
0\to I^q\to A\to A/I^q\to 0
$$
induces an exact sequence
$$
0\to\Hom_A(A/I^q,R)\tto\al\Hom_A(A,R)\tto\beta\Hom_A(I^q,R)\tto\ga\Ext_A^1(A/I^q,R)\to 0.
$$
Now $\al$ is an isomorphism, therefore $\beta$ is zero and $\ga $ is an isomorphism again.
We get an isomorphism
$$
\Hom_R(I^q/I^{q+1},R)\ \cong\ \Hom_A(I^q,R)\ \cong\ \H_q^1(\Ga,R).
$$
\mathqed

\begin{proposition}\label{Prop2.4}
Suppose that $R$ is a field.
Let $\Ga$ be finitely generated and let $V$ be an $R[\Ga]$-module with
$$
H^1(\Ga,V)\= 0.
$$
Then for every $q\ge 1$ there is a natural isomorphism
$$
\ol\H_{q+1}^0(\Ga,V)\ \tto\cong\  
\H_{q}^1(\Ga,R)\otimes V^\Ga.
$$ 
\end{proposition}

\prf
By Proposition \ref{1.1} we get an exact sequence
$$
0\to \H_{q}^0(\Ga,V)\to \H_{q+1}^0(\Ga,V)\to \Hom_R(I^q/I^{q+1},\H^0(\Ga,V))\to 0.
$$
Since $\Ga$ is finitely generated, the free $R$-module $I^q/I^{q+1}$ is finite-dimensional, hence we get
$$
\Hom_R(I^q/I^{q+1},H^0(\Ga,V))\ \cong\ \Hom_R(I^q/I^{q+1},R)\otimes H^0(\Ga,V).
$$
The claim now follows from Lemma \ref{Prop2.1.3}.
\qed

\begin{lemma}\label{Lem2.5}
For $q\ge 1$ the natural map
$$
\H_{q}^1(\Ga,R)\ \to\ \H_{q+1}^1(\Ga,R)
$$
is the zero map.
\end{lemma}

\prf
In Lemma \ref{1.1} we put $V=R$ and get the exact sequence
$$
\Hom_R(I^q/I^{q+1},R)\tto\al\H_{q}^1(\Ga,R)\to H_{q+1}^1(\Ga,R).
$$
We claim that $\al$ is an isomorphism.
Note that $\al$ is the restriction of the connection morphism attached to the exact sequence
$$
0\to I^q/I^{q+1}\to A/I^{q+1}\to A/I^q\to 0,
$$
which is related to the sequence
$$
0\to I^q\to A\to A/I^q\to 0.
$$
The connection morphism $\tilde\al$ of the latter was shown in Lemma \ref{Prop2.1.3} to be an isomorphism.
Taking residue classes modulo $I^{q+1}$ maps the first exact sequence to the second.
This results in a map of the corresponding long exact $\Ext$-sequences, a part of which gives the commutative diagram
$$
\begin{diagram}
\node{\Hom_A(I^q/I^{q+1},R)}
		\arrow{e,t}\al
		\arrow{s,r}{\cong}
	\node{\H_{q}^1(\Ga,R)}
		\arrow{s,r}\cong\\
\node{\Hom_A(I^q,R)}
		\arrow{e,t}\cong
	\node{\H_{q}^1(\Ga,R).}
\end{diagram}
$$
It follows that $\al$ is an isomorphism, indeed.
\qed

\section{Hecke pairs and smooth modules}
A \e{Hecke pair} is a pair $(\CG,\Ga)$ of a group $\CG$ and a subgroup $\Ga$ such that for every $g\in \CG$ the set $\Ga g\Ga/\Ga$ is finite.
We also say that $\Ga$ is a \e{Hecke subgroup} of $\CG$.

Two subgroups $\Ga,\La$ of a group $H$ are called \e{commensurable}, written $\Ga\sim\La$,  if the intersection $\Ga\cap\La$ has finite index in both.
Commensurability is an equivalence relation which is preserved by automorphisms of $H$.

The \e{commensurator} of a group $\Ga\subset H$ is
$$
\comm(\Ga)\df \{ h\in H:\Ga \text{ and } h\Ga h^{-1}\text{ are commensurable}\}.
$$

\begin{lemma}
The commensurator $G=\comm(\Ga)$ is a subgroup of $H$.
It is the largest subgroup such that $(G,\Ga)$ is a Hecke pair.
More precisely,
$$
\comm(\Ga)\=\{ h\in H: |\Ga h\Ga/\Ga|,|\Ga\bs\Ga h\Ga|<\infty\}.
$$
\end{lemma}

\prf
Let $G=\comm(\Ga)$ and $g\in G$.
Then $\Ga\sim g\Ga g^{-1}$.
By conjugating we get $g^{-1}\Ga g\sim\Ga$, hence $g^{-1}\in G$.
Next let $h\in G$ as well.
Then $\Ga\sim h\Ga h^{-1}$ and hence $g\Ga g^{-1}\sim gh\Ga (gh)^{-1}$ so that $\ga\sim g\Ga g^{-1}$ implies $\Ga\sim gh\Ga(gh)^{-1}$, which means that $gh\in G$, so $G$ is indeed a subgroup.

To see the identity claimed in the lemma, one simply observes that for every $h\in H$ the natural map
$
\Ga/\Ga\cap h\Ga h^{-1}\ \to\ \Ga h\Ga/\Ga,
$
mapping the class of $\ga$ to the class of $\ga h$ is a bijection.
\qed

Let $\CG$ be a group.
By a $\CG$-module we shall henceforth mean an $R[\CG]$-module.
If $\CG$ is a totally disconnected topological group, an element $v$ of a $\CG$-module $V$ is called \e{smooth} if it is stabilized by some open subgroup of the topological group $\CG$.
The set $V^\infty$ of all smooth elements is a submodule and the module $V$ is called \e{smooth} if $V=V^\infty$.

Drop the condition that $\CG$ be a topological group and let $(\CG,\Ga)$ be a Hecke-pair.
A \e{congruence subgroup} of $\Ga$ is any subgroup which contains a group of the form
$$
\Ga\cap g_1\Ga g_1^{-1}\cap\dots\cap g_n\Ga g_n^{-1}
$$
for some $g_1,\dots, g_n\in \CG$.
As $(\CG,\Ga)$ is a Hecke pair, every congruence subgroup has finite index in $\Ga$.
Note that this definition of a congruence subgroup coincides with the one given in \cite{HOAut}.
For every congruence subgroup $\Sigma$ equip the set $\CG/\Sigma$ with the discrete topology and consider the topological space
$$
\what\CG\df \lim_{\stack\leftarrow\Sigma}\CG/\Sigma,
$$
where the limit is taken over all congruence subgroups $\Sigma$.

\begin{lemma}
\begin{enumerate}[\rm (a)]
\item The intersection of all congruence subgroups $N=\bigcap_\Sigma\Sigma$ is a normal subgroup of $\CG$.
\item The natural map $p:\CG\to\what\CG$ factors through the injection $\CG/N\hookrightarrow \what \CG$ and has dense image.
\item The group multiplication on $\CG/N$ extends by continuity to $\what \CG$ and makes $\what\CG$ a totally disconnected locally compact group.
\end{enumerate}
\end{lemma}

We call $\what\CG$ the \e{congruence completion} of $\CG$.
Although the notation doesn't reflect this, the completion $\what\CG$ depends on the choice of the Hecke subgroup $\Ga$.
A Hecke subgroup $\Ga$ is called \e{effective}, if the normal subgroup $N$ above is trivial.

\prf
(a) Let $n\in N$ and let $g\in \CG$.
For a given congruence subgroup $\Sigma$ we have that $n\in\Sigma\cap g^{-1}\Sigma g$, and so $gng^{-1}\in\Sigma$.
As $\Sigma$ varies, we find $gng^{-1}\in N$.

(b) Let $g,g'\in\CG$ with $p(g)=p(g')$.
This means that $g\Sigma=g'\Sigma$ for every congruence subgroup and so $gN=g'N$. 
For given $(g_\Sigma)_\Sigma\in\what\CG$ the sets $U_\Sigma=\{ h\in\what\CG: h_\Sigma\Sigma=g_\Sigma\Sigma\}$ form a neighborhood base.
Clearly the element $g_\Sigma\in\CG$ is mapped into $U_\Sigma$, so the image of $p$ is dense.

(c) Let $\bar g=(g_\Sigma)_\Sigma\in\what\CG$.
Then the net $(p(g_\Sigma))_\Sigma$ converges to $\bar g$.
For $\bar h=(h_\Sigma)\in\what\CG$ it is easy to see that the net $(p(g_\Sigma h_\Sigma))_\Sigma$ converges in $\bar\CG$.
We call the limit $\bar g\bar h$.
This multiplication has the desired properties.
\qed

The initial topology defined by $p$ on $\CG$ makes $\CG$ a topological group with the congruence groups forming a unit neighborhood base.
Clearly every smooth $\bar\CG$-module is a smooth $\CG$-module by restriction.
But also the converse is true: Every smooth $\CG$-module extends uniquely to a smooth $\bar\CG$-module and these two operations of restriction and extension are inverse to each other.
Let
$$
\what\H_{q}^p(V)\df \lim_{\stack\to\Sigma}\H_{q}^p(\Sigma,V),
$$
where the limit is taken over all congruence subgroups of $\Ga$. For $p=0$ we also define 
$$
\what{\ol\H}_{q}^0(V)=\lim_{\stack{\to}{\Sigma}}\ol\H_q^0(\Sigma,V).
$$
Note that $\what{\ol\H}_{0}^0(V)=V^\infty$.
For $g\in\CG$, the map induced by $g$:
$$
 \H_{q}^p(\Sigma,V)\to  \H_{q}^p(g\Sigma g^{-1},V)\tto\res \H_{q}^p(\Sigma\cap g\Sigma g^{-1},V)
$$
defines an action of $\CG$ on $\what\H_{q}^p(V)$.

Assume from now on that $\Ga$ is finitely generated and $R$ is a field.
Then every finite-index subgroup $\Sigma\subset\Ga$ is finitely generated as well, so Lemma \ref{1.1*} applies to all modules $V$.

Let $R$ be a field and consider the order-lowering map
$$
\ol \H_{q}^0(\Sigma,V)\hookrightarrow\Hom(\Sigma,R)\otimes\ol \H_{q-1}^0(\Sigma,V).
$$
By iteration we get an injection
$$
\ol \H_{q}^0(\Sigma,V)\hookrightarrow\Hom(\Sigma,R)^{\otimes (q-1)}\otimes V^\Sigma.
$$
Taking limits we get an injection
$$
\what{\ol\H}_{q}^0(V)\hookrightarrow \(\lim_{\stack\to\Sigma}\Hom(\Sigma,R)\)^{\otimes (q-1)}\otimes V^\infty.
$$
Write
$$
\what\Hom(\Ga,R)\=\lim_{\stack\to\Sigma}\Hom(\Sigma,R).
$$
We have shown:

\begin{proposition}\label{4.2}
If $\Ga$ is finitely generated and $R$ is a field,
then there is a natural injection of smooth modules 
$$
\what{\ol\H}_{q}^0(V)\ \hookrightarrow\ 
\what\Hom(\Ga,R)^{\otimes (q-1)}\otimes V^\infty.
$$ 
\end{proposition}

There are examples when this map is not surjective.

\begin{proposition}\label{Prop3.4}
Let $(\CG,\Ga)$ be a Hecke-pair.
Assume $\Ga$ to be finitely generated and let $V$ be an $R[\Ga]$-module such that
$$
\H^1(\Sigma,V)=0
$$
for every congruence subgroup $\Sigma$ of $\Ga$.
Then for every $q\ge 1$ there is an isomorphism of $\CG$-modules
$$
\what{\ol\H}_{q}^0(V)\tto\cong \what\H_{q-1}^1(R)\otimes V^\infty.
$$
\end{proposition}

\prf
This is a consequence of Proposition \ref{Prop2.4}.
\qed

\section{Higher order $L^2$-spaces}
Let $G$ be a locally compact group and $\Ga$ a countable discrete subgroup.
We consider the action of $\Ga$ on $G$ by left translations.

\begin{lemma}
The group $\Ga$ is closed in $G$.
It acts \e{strongly discontinuously} in the sense that for every compact subset $C\subset G$ the set 
$$
\{\ga\in\Ga: C\cap\ga C\ne \emptyset\}
$$
is finite.
\end{lemma}

\prf
It is known that discrete subgroups are closed \cite{HA2}.
Further, discreteness of $\Ga$ implies that $\Ga$ meets every compact set in only finitely many points.
With $C$ the set $CC^{-1}$ is also compact, which gives the claim.
\qed

A \e{fundamental set} is a set $F\subset G$ of representatives for the quotient $\Ga\bs G$ which is measurable.

A fundamental set $F$ is called \e{locally finite} if there exists an open neighborhood $U$ of the closure $\ol F$ such that $U$ only meets finitely many translates of $\ol F$, i.e., the set
$$
\{ \ga\in\Ga: \ga \ol F\cap U\ne \emptyset\}
$$
is finite.

\begin{lemma}\label{Lem4.1}
Any measurable set $M\subset G$ which contains a set of representatives for $\Ga\bs G$, contains a fundamental set.
\end{lemma}

\prf
According to a theorem of Feldman and Greenleaf \cite{FG} there exists a fundamental set $\tilde F\subset G$.
We now enumerate the elements of $\Ga$, so
$$
\Ga\=\{1=\ga_1,\ga_2,\ga_3,\dots\}.
$$
We set $F_1=M\cap \tilde F$ and iteratively we define
$$
F_{j+1}\= F_j\cup \left[ \left( M\cap\ga_{j+1} \tilde F\right)\sm \Ga M_j\right].
$$
Then $F=\bigcup_{j=1}^\infty F_j$ is a fundamental set contained in $M$.
\qed

Two fundamental sets $F_1,F_2$ are called \e{compatible}, if there is a finite subset $E\subset\Ga$ such that
$$
F_1\subset  E F_2\quad\text{and}\quad F_2\subset  E F_1.
$$
Let $M$ be the space of all measurable functions modulo nullfunctions on $G$ and let $M_q$ be the subspace of all $f\in F$ such that $I^{q}f=0$, where $I$ is the augmentation ideal in $\C[\Ga]$.

\begin{proposition}
Fix a fundamental set $F\subset G$.
Assume that $\Ga$ is finitely generated and 
let $S$ be a finite symmetric set of generators of $\Ga$ which is supposed to contain the unit element.
Then any $f\in M_{q}$ is uniquely determined by its restriction to
$$
S^{q-1}F\=\bigcup_{s_1,\dots,s_{q-1}\in S}s_1\cdots s_{q-1} F.
$$
here we interpret $S^0F$ as $F$ itself.
The space 
$$
M_{q}\cap L^2(S^{q-1}F)
$$
does neither depend on the choice of $S$ nor on the choice of $F$ in a given compatibility class.
We denote this space by $L^2_{q}(F)$ or, if there is a given choice of $F$, by $L^2_{q}(\Ga\bs X)$.

Also, the topology on $L_{q}^2(F)$, given by the $L^2$-structure, is independent of the choices in the same way as the space itself.
\end{proposition}

\prf
We have to show that any $f\in M_q$ which vanishes on $S^{q-1}F$, is zero.
We use induction on $q$.
The case $q=1$ is clear.
Let $q\ge 2$ and set $\bar M_q=M_q/M_{q-1}$.
Consider the \e{order lowering operator}
$$
\La:M_q\to \Hom(\Ga,\bar M_{q-1})\cong\Hom(\Ga,\C)\otimes\bar M_{q-1}
$$
given by
$$
\La(f)(\ga)\=(\ga-1) f.
$$
The kernel of $\La$ is $M_{q-1}$.
Now assume $f\in M_q$ vanishes on $S^{q-1}F$.
For every $s\in S$ we have $(s-1)f|_{S^{q-1}F}=0$ and hence, by induction hypothesis, $(s-1)f=0$.
As $S$ generates $\Ga$, this means $\La(f)=0$ and so $f\in M_{q-1}$, so again by induction hypothesis we get $f=0$.

We next show
$$
M_q\cap L^2(S^qF)=M_q\cap L^2(S^{q+j}F)
$$
for every $j\ge 0$.
The inclusion ``$\supset$'' is clear.
We show ``$\subset$'' by induction on $q$.
For $q=1$ the claim is clear.
So assume $q\ge 2$ and the claim proven for $q-1$.
We show $M_q\cap L^2(S^{q-1+j}F)\subset M_q\cap L^2(S^{q+j}F)$ for every $j\ge 0$.
For this let $f\in M_q\cap L^2(S^{q-1+j}F)$ and let $s\in S$.
Then $f(s^{-1} z)=f(z)+f(s^{-1} z)-f(z)=f(z)+(s-1)f(z)$.
The function $f(z)$ is in $L^2(S^{q-1+j}F)$ and the function $(s-1)f(z)$ is in $M_{q-1}\cap L^2(S^{q-2+j}F)$ the latter space equals $M_{q-1}\cap L^2(S^{q-1+j}F)$ by induction hypothesis.
It follows that $f(s^{-1} z)$ is in $L^2(S^{q-1+j}F)$, so $f\in L^2(sS^{q-1+j}F)$.
Since this holds for every $s\in S$, we get $f\in L^2(S^{q+j}F)$ as claimed.

From this we conclude the independence of $S$.
Let $S_1$ be a second finite symmetric set of generators containing the unit.
Then there exists $j\ge 0$ such that $S_1^q\subset S^{q+j}$.
Hence it follows that $M_q\cap L^2(S_1^{q-1}F)\subset M_q\cap L^2(S^{q-1+j}F)=M_q\cap L^2(S^{q-1}F)$.
By symmetry we get $M_q\cap L^2(S_1^{q-1}F)=M_q\cap L^2(S^{q-1}F)$ as claimed.
We have to show independence of the choice of $F$.
So let $F'$ be another fundamental set which is compatible with $F$.
Then there exists $j\ge 0$ such that $S^qF'\subset S^{q+j}F$ and we prove independence as above.

Finally, for the topology, let $(f_\nu)$ be a sequence in $M_q\cap L^2(S^{q-1}F)$ tending to zero.
We argue as in the proof of $M_q\cap L^2(S^{q-1}F)=M_q\cap L^2(S^{q-1+j}F)$ to show that $(f_\nu)$ also tends to zero in $L^2(S^{q-1+j}F)$.
From here on the proof of independence of choices for the topology on the space $L^2_{q}(F)$ proceeds as above.
\qed

We will remove the dependence on the choice of a fundamental set by giving a canonical choice (up to compatibility) in the case when $\Ga\bs G$ is compact.

\begin{lemma}\label{Lem4.5}
There exists a relatively compact fundamental set if and only if $\Ga\bs G$ is compact.
Each relatively compact fundamental set is locally finite and any two relatively compact fundamental sets are compatible.
\end{lemma}

\prf
Let $\pi:G\to\Ga\bs G$ be the projection.
This map is continuous and open.
Let $F$ be a relatively compact fundamental set, then $\Ga\bs G$ is the image under $\pi$ of the compact set $\ol F$, hence $\Ga\bs G$ is compact.

For the converse assume that $\Ga\bs G$ is compact.
For each $x\in G$ fix a relatively compact open neighborhood $U_x$ of $x$.
The sets $\pi(U_x)$ form an open covering of $\Ga\bs G$, so finitely many suffice.
This means there are $x_1,\dots, x_n\in G$ such that the compact set $\tilde F=\ol U_{x_1}\cup\dots\cup \ol U_{x_n}$ contains a set of representatives.
By Lemma \ref{Lem4.1} it contains a fundamental set which is necessarily relatively compact.

For the second assertion let $F$ be a relatively compact fundamental set.
Let $U\supset F$ be a relatively compact open neighborhood of $\ol F$.
As $\Ga$ acts strongly discontinuously, it follows that $\ol U$ meets only finitely many translates of $\ol F$ and hence $U$ only meets finitely many translates of $\ol F$.

Finally, let $F_1,F_2$ be two relatively compact fundamental sets.
Let $U$ be an open neighborhood of $F_1$ that only meets finitely many translates of $F_1$.
Then the family $(\ga U)_{\ga\in\Ga}$ is an open covering of $G$, so there are $\ga_1,\dots\ga_n$ such that $\ol F_2\subset \ga_1 U\cup\dots\cup\ga_n U$.
Now $U$ is contained in a finite number of translates of $F_1$, so there exists a finite set $E\subset \Ga$ such that $F_2\subset EF_1$.
By symmetry, we get the other direction, too.
\qed

From now on we assume that $\Ga$ is cocompact, i.e. $\Ga\bs G$ is compact.
In this case we say that $\Ga$ is a \e{uniform lattice} in $G$.
We consider the space
$
L^2_{\loc}(G)
$
of all locally square-integrable functions on $G$ (modulo nullfunctions).
So we have $f\in L^2_\loc(G)$ if for every $x\in G$ there exists an open neighborhood $U$ such that $f|_U\in L^2(U)$.
This is equivalent to saying that $f|_K\in L^2(K)$ for every compact subset $K\subset G$.
We have
$$
L_{q}^2(\Ga\bs G)\= H_{q}^0(\Ga,L^2_{\loc}(G)).
$$

\begin{definition}
The elements of $L_q^2(\Ga\bs G)$ will be called \e{automorphic forms of order $q$}.
Automorphic forms of order 1 are classical automorphic forms. Those of order $\ge 2$ are also referred to as \e{automorphic forms of higher order}.
\end{definition}

\begin{proposition}\label{Prop4.3}
For every $q\ge 0$ we have
$$
H_{q}^1(\Ga,L^2_{\loc}(G))\= 0.
$$
\end{proposition}

\prf
Let $F$ be a relatively compact fundamental set.
By Lemma \ref{1.1} (b)  it suffices to consider the case $q=0$.
Let $\al:\Ga\to L_{\loc}^2(G)$ be a 1-cocycle, i.e. a map that satisfies $\al(\ga\tau)=\ga\al(\tau)+\al(\ga)$.
We set
$$
f(x)\=\sum_{\tau\in\Ga}\al(\tau)(x)\1_F(\tau^{-1}x).
$$
A simple computation shows $(1-\ga)f=\al(\ga)$ and since $\al\in L^2_{\loc}(G)$ it follows that $f\in L^2_{\loc}(G)$.
The proposition follows.
\qed

\begin{proposition}\label{Prop4.6}
For every $q\ge 1$ there is a natural exact sequence of continuous linear maps,
$$
0\to L^2_{q}(\Ga\bs G)\to  L^2_{q+1}(\Ga\bs G)
\to L^2(\Ga\bs G)^{N_{\Ga}(q)}\to 0.
$$
\end{proposition}

\prf
This follows from Lemma \ref{1.1} together with Proposition \ref{Prop4.3}.
\qed

\section{Lie groups}
In this section $G$ will be a Lie group.
In this case we write
$$
C_q^\infty(\Ga\bs G)\= H_q^0(\Ga,C^\infty(G)).
$$

\begin{proposition}
For a Lie group $G$ and a discrete subgroup $\Ga$ we have
$$
H^1(\Ga, C^\infty(G))\= 0.
$$
Consequently, for every $q\ge 1$ there is a natural exact sequence
$$
0\to C_{q}^\infty(\Ga\bs G)\to C_{q+1}^\infty(\Ga\bs G)\to C^\infty(\Ga\bs G)^{N_\Ga(q)}\to 0.
$$
\end{proposition}

\prf
This is a a part of Proposition 2.3.2 in \cite{HOinv}.
\qed

Let $\pi\in\hat G$.
A representation $(\beta,V_\beta)$ of $G$ is said to be \e{of type $\pi$}, if it is of finite length and every irreducible subquotient is isomorphic to $\pi$.
For a representation $(\eta,V_\eta)$ we define the \e{$\pi$-isotype} as
$$
V_\eta(\pi)\df \ol{\sum_{\stack{V_\beta\subset V_\eta}{\beta{\rm\ of\ type\ }\pi}}V_\beta},
$$
where the sum runs over all subrepresentations $V_\beta$ of type $\pi$.
If $\pi\ne \pi'$, then $V_\eta(\pi)\cap V_\eta(\pi')=0$.

If $\Ga$ is a uniform lattice, then the representation of $G$ on $L^2(\Ga\bs G)$ is a direct sum of irreducible representations, each occuring with finite multiplicity, i.e.,
$$
L^2(\Ga\bs G)\ \cong\ \bigoplus_{\pi\in\what G}m_\Ga(\pi)\pi,
$$
with $m_\Ga(\pi)\in\N_0$.

\begin{theorem}[Spectral decomposition]\label{thm5.2}
Let $G$ be a semisimple Lie group and $\Ga$ a uniform lattice in $G$.
We write $V_q=\H_q^0(\Ga,\Omega)$, where $\Omega=L_\loc^2(G)$ or $\Omega=C^\infty(G)$.
For every $q\ge 1$ there is an isotypical decomposition
$$
V_q\=\ol{\bigoplus_{\pi\in\hat G_\Ga}}V_q(\pi),
$$
and each $V_q(\pi)$ is of type $\pi$ itself.
The exact sequence of Proposition \ref{Prop4.6} induces an exact sequence
$$
0\to V_{q}(\pi)\to V_{q+1}(\pi)\to V_1(\pi)^{m_\Ga(\pi)N_{\Ga}(q)}\to 0
$$
for every $\pi\in\hat G$.
\end{theorem}

\prf
For $\Omega= C^\infty(G)$ this is Theorem 2.3.5 of \cite{HOinv}.
For $\Omega= L^2_\loc(G)$ the result follows from the latter by density arguments.
\qed

\section{Sheaf cohomology}
Let $Y$ be a topological space which is path-connected and locally simply connected.
Let $\Ga$ be the fundamental group of $Y$ and let $X\tto\pi Y$ be the universal covering.

For a sheaf $\CF$ on $Y$ define
$$
H_{q}^0(Y,\CF)\df H_{q}^0(\Ga,H^0(X,\pi^*\CF)).
$$
Let $\Mod(R)$ be the category of $R$-modules, let $\Mod_R(Y)$ be the category of sheaves of $R$-modules on $Y$, and let $\Mod_R(X)_\Ga$ be the category of sheaves over $X$ with an equivariant $\Ga$-action.
Then $\H_{q}^0(Y,\cdot)$ is a left exact functor from $\Mod_R(Y)$ to $\Mod(R)$.
We denote its right derived functors by $\H_{q}^p(Y,\cdot)$ for $p\ge 0$.

\begin{proposition}\label{Prop8.1} Assume that the universal cover $X$ is contractible.
\begin{enumerate}[\rm (a)]
\item For each $p\ge 0$ one has a natural isomorphism $\H_{1}^{p}(Y,\CF)\ \cong\ \H^p(Y,\CF)$.
\item If a sheaf $\CF$ is $\H^0(Y,\cdot)$-acyclic, then it is $\H_{q}^0(Y,\cdot)$-acyclic for every $q\ge 1$.
\end{enumerate}
\end{proposition}

Note that part (b) allows one to use flabby or fine resolutions to compute higher order cohomology.

\prf
We decompose the functor $\H^0_{q}(Y,\cdot)$ as a composition of functors
$$
\Mod_R(Y)\tto{\pi^*}\Mod_R(X)_\Ga\tto{\H^0(X,\cdot)}\Mod(R[\Ga])\tto{\H^0_{q}(\Ga,\cdot)}\Mod(R).
$$
The functor $\pi^*$ is exact and maps injectives to injectives.
We claim that $\H^0(X,\cdot)$ has the same properties.
For the exactness, consider the commutative diagram
$$
\begin{diagram}\divide\dgARROWLENGTH by2
\node{\Mod_R(X)_\Ga}\arrow{e,t}{\H^0}\arrow{s,r}f
	\node{\Mod(R[\Ga])}\arrow{s,r}f\\
\node{\Mod_R(X)}\arrow{e,t}{\H^0}
	\node{\Mod(R),}
\end{diagram}
$$
where the vertical arrows are the forgetful functors.
As $X$ is contractible, the functor $\H^0$ below is exact.
The forgetful functors have the property, that a sequence upstairs is exact if and only if its image downstairs is exact.
This implies that the above $\H^0$ is exact.
It remains to show that $\H^0$ maps injective objects to injective objects.
Let $\CJ\in\Mod_R(X)_\Ga$ be injective and consider a diagram with exact row in $\Mod(R[\Ga])$, 
$$
\begin{diagram}\divide\dgARROWLENGTH by2
\node{0}\arrow{e}
	\node{M}\arrow{e}\arrow{s,r}{\ph}
	\node{N}\\
\node[2]{\H^0(X,\CJ).}
\end{diagram}
$$
The morphism $\ph$ gives rise to a morphism $\phi:M\times X\to \CJ$, where $M\times X$ stands for the constant sheaf with stalk $M$.
Note that $\H^0(X,\phi)=\ph$.
As $\CJ$ is injective, there exists a morphism $\psi:N\times X\to \CJ$ making the diagram
$$
\begin{diagram}\divide\dgARROWLENGTH by2
\node{0}\arrow{e}
	\node{M\times X}\arrow{e}\arrow{s,r}{\phi}
	\node{N\times X}\arrow{sw,r}{\psi}\\
\node[2]{\CJ}
\end{diagram}
$$
commutative.
This diagram induces a corresponding diagram on the global sections, which implies that $\H^0(X,\CJ)$ is indeed injective.

For a sheaf $\CF$ on $Y$ it follows that
$$
\H^p(Y,\CF)\=R^p(\H^0(Y,\CF))\=R^p\H^0_{q}(\Ga,\CF)\circ \H^0_\Ga\circ\pi^*\= \H_{q}^p(Y,\CF).
$$
Now let $\CF$ be acyclic.
Then we conclude $\H_{0}^p(Y,\CF)=0$ for every $p\ge 1$, so the $\Ga$-module $V=\H^0(X,\pi^*\CF)$ is $\Ga$-acyclic.
The claim follows from Lemma \ref{1.1}.
\qed

Let $Z$ be a closed subset of $Y$ and let $U$ be its complement.
Let $i:Z\hookrightarrow X$ and $j:U\hookrightarrow X$ denote the inclusions.

\begin{proposition}
Let $\CF$ be a sheaf on $Y$.
The exact sequence
$$
0\to j_!(\CF|_U)\to \CF\to i_*(\CF|_Z)\to 0
$$
gives rise to a long exact sequence of cohomology groups:
\begin{multline*}
0\to \H_{q,c}^0(U,\CF)\to \H_{q}^0(Y,\CF)\to\H_{q}^0(Z,\CF)\\
\to \H_{q,c}^1(U,\CF)\to \H_{q}^1(Y,\CF)\to\H_{q}^1(Z,\CF)\\
\dots \to \H_{q,c}^p(U,\CF)\to \H_{q}^p(Y,\CF)\to\H_{q}^p(Z,\CF)\to\dots
\end{multline*}
Here $\H_{q,c}^0(U,\CF)$ stands for the group of all sections in $\H_{q}^0(Y,\CF)$ whose support is contained in $\pi^{-1}(U)$ and $\H_{q}^p(Y,\CF)$ for $p\ge 1$ is the corresponding derived functor.
\end{proposition}

\prf
It is easy to see that $\H_{q}^0(Y,j_!(\CF|_U))\=\H_{q,c}^0(U,\CF)$. 
The functor $\CF\mapsto j_!(\CF|_U)$ is exact and maps products of injective skyscraper sheaves to products of injective skyscraper sheaves.
As any sheaf has a resolution consisting of products of injective skyscraper sheaves (which are injective themselves), it follows that $\H_{q}^p(Y,j_!(\CF|_U))\=\H_{q,c}^p(U,\CF)$ holds for every $p\ge 0$.
This proves the assertion.
\qed

Assume that $X$ is contractible.
Then $Y$ is a classifying space for the group $\Gamma$.
So any $R[\Ga]$-module $V$ gives rise to a locally constant sheaf $\CV$ of $R$-modules over $Y$ such that $H^p(\Ga,V)\cong H^p(Y,\CV)$ holds for every $p\ge 0$.

\begin{proposition}
We have natural isomorphisms
$$
\H_q^p(\Ga,V)\ \cong\ \H_q^p(Y,\CV)
$$
for all $p\ge 0, q\ge 1$.
\end{proposition}

\prf
Fix $q\ge 1$.
The functors $V\mapsto \H_q^p(\Ga,V)$ form a universal $\delta$-functor on $\Mod(R[\Ga])$.
Let $Sh(Y)$ be the category of sheaves of $R$-modules over $Y$, so $\CF\mapsto \H_q^p(Y,\CF)$ is a universal $\delta$-functor on $Sh(Y)$.
The functor of sheafification $\Sheaf:\Mod(R[\Ga])\to Sh(Y)$, which to a module  $V$ attaches the locally constant sheaf 
$\CV=\Sheaf(V)$, is exact.
So $V\mapsto \H_q^p(Y,\Sheaf(V))$ is a  $\delta$-functor on $\Mod(R[\Ga])$.
We have to show universality, which we do as usual by showing erasability of $\H^p$ for $p\ge 1$.
For $V\in\Mod(R[\Ga])$ let
$$
IV=\{ \al:\Ga\to V\}
$$
the module of all maps from $\Ga$ to $M$.
This is a  $\Ga$-module via
$$
\ga. \al(\tau)\= \ga(\al(\ga^{-1}\tau)).
$$
Mapping $v\in V$ to the constant map $v$, one gets an embedding $V\hookrightarrow IV$.
By Proposition \ref{Prop8.1} we have to show
$$
\H^p(Y,\Sheaf(IV))=0
$$ 
for $p\ge 1$.
Let $\pi: X\to Y$ be the projection.
Then
$$
\Sheaf(IV)\=\pi_*\tilde V,
$$
where $\tilde V$ stands for the constant sheaf $V$ on $X$.
We now show, that $\pi_*\tilde V$ is acyclic.
To this end we resolve $\tilde V$ with special acyclic sheaves, which are products of skyscraper sheaves with injective stalks.
The images under  $\pi_*$ of these are again products of skyscraper sheaves with injective stalks, hence acyclic.
The global sections above and below are the same, so the cohomologies agree, i.e., we have
$$
\H^p(Y,\pi_*\tilde V)\= \H^p(X,\tilde V).
$$
As $X$ is contractible, the right hand side is zero.
\qed

\section{Cohomology of lattices in Lie groups}
In this section we set $R=\C$.
Let $G$ be a semisimple Lie group with compact center and finitely many connected components.
Let $\Ga\subset G$ be a uniform lattice which is torsion-free.
Let $Y=\Ga\bs X$, then $\Ga$ is the fundamental group of the manifold $Y$, and the universal covering $X$ of $Y$ is contractible.
This means that we can apply the results of the last section.
Compare the next result to the classical case \cite{Schw}.

\begin{theorem}
Let $(\sigma,E)$ be a finite dimensional irreducible representation of $G$.
There is a natural isomorphism
\begin{eqnarray*}
\H_{q}^p(\Ga,E)
&\cong& \H^p(\g,K,\H_{q}^0(\Ga,C^\infty(G))\otimes E)\\
&=& \bigoplus_{\stack{\pi\in\hat G}{\chi_\pi=\chi_{E^*}}}
\H^p(\g,K,C_q^\infty(\Ga\bs G)(\pi)\otimes E),
\end{eqnarray*}
where $\H(\g,K,.)$ denotes $(\g,K)$-cohomology and the direct sum is finite.
\end{theorem}

\prf
Let $\CF_E$ be the locally constant sheaf on $Y$ corresponding to $E$.
Let $\Omega_Y^p$ be the sheaf of complex valued $p$-differential forms on $Y$.
Then $\Omega_Y^p\otimes\CF_E$ is the sheaf of $\CF_E$-valued differential forms.
These form a fine resolution of $\CF_E$:
$$
0\to\CF_E\to\C^\infty\otimes\CF_E\tto{d\otimes 1}\Omega_Y^1\otimes\CF_E\to\dots
$$
Since $\pi^*\Omega_Y^\bullet=\Omega_X^\bullet$, 
we conclude that $\H^p_{q}(\Ga,E)$ is the cohomology of the complex $\H_{q}^0(\Ga,\H^0(X,\Omega_X^\bullet\otimes E))$.
Let $\g$ and $\k$ be the Lie algebras of $G$ and $K$ respectively, and let $\g=\k\oplus\p$ be the Cartan decomposition.
Then $\H^0(X,\Omega^p\otimes\CF_E)\=(C^\infty(G)\otimes\bigwedge^p\p)^K\otimes E$. 
Mapping a form $\omega$ in this space to $(1\otimes x^{-1})\omega(x)$ one gets an isomorphism to $(C^\infty(G)\otimes\bigwedge^p\p\otimes E)^K$, where $K$ acts diagonally on all factors and $\Ga$ now acts on $C^\infty(G)$ alone.
The first claim follows.
For the second, we replace $C_q^\infty(\Ga\bs G)$ by the sum of its isotypes according to Theorem \ref{thm5.2}.
The direct sum can be pulled out of cohomology 
to get
$$
\H_q^p(\Ga,E)\ \cong\ \bigoplus_{\pi\in\what G}\H(\g,K,V_q(\pi)\otimes E),
$$
where $V_q=C_q^\infty(\Ga\bs G)$.
For each $\pi$ there is a $G$-stable finite filtration
$$
0=F_0\subset\dots\subset F_n-V_q(\pi)
$$
such that $F_j/F_{j-1}\cong\pi^\infty$.
The exact sequence
$$
0\to F_{j-1}\to F_j\to \pi^\infty\to 0
$$
induces a long  exact cohomology sequence
\begin{multline*}
0\to \H^0(\g,K,F_{j-1}\otimes E)\to\H^0(\g,K,F_j\otimes E)\to \H^0(\g,K,\pi^\infty\otimes E)\to\\
\to\H^1(\g,K,F_{j-1}\otimes E)\to\H^1(\g,K,F_j\otimes E)\to \H^1(\g,K,\pi^\infty\otimes E)\to\\
\dots\to\H^p(\g,K,F_{j-1}\otimes E)\to\H^p(\g,K,F_j\otimes E)\to \H^p(\g,K,\pi^\infty\otimes E)\to\dots
\end{multline*}
Now assume $\chi_\pi\ne\chi_{E^*}$.
By Theorem 4.1 of \cite{BorelWallach} 
we conclude $\H^p(\g,K,\pi^\infty\otimes E)=0$ for all $p$.
By induction on $j$ one gets $\H^p(\g,K,F_j\otimes E)=0$ for all $p$ and all $j$.
Hence in the above sum only those finitely many summands with $\chi_\pi=\chi_{E^*}$ remain.
\qed

Taking $q$-th powers gives a natural surjective map
$$
\(I/I^2\)^{\otimes q}\to I^q/I^{q+1}.
$$
This dualizes to an injection
$$
\Hom_\C(I^q/I^{q+1},\C)\ \hookrightarrow\ \Hom_\C(I/I^2,\C)^{\otimes q}.
$$

For $p=1$ and $q\ge 1$ we thus get,
\begin{eqnarray*}
\H_{q-1}^1(\Ga,\C) &\cong& \Hom_\C(I^q/I^{q+1},\C)\\
&\hookrightarrow& \Hom_\C(I/I^2,\C)^{\otimes q}\\
&\cong& \H^1(\Ga,\C)^{\otimes q}\\
&\cong& \(\bigoplus_{\stack{\pi\in\what G}{\chi_\pi=\chi_\C}}m_\Ga(\pi)\H^1(\g,K,\pi^\infty)\)^{\otimes q}
\end{eqnarray*}

\section{Arithmetic groups}
In this section we put $R=\C$.
Let $\G$ be a linear algebraic group over $\Q$ which is simple and simply connected and such that $G=\G(\R)$ has no compact component.
By strong approximation, the group $\CG=\G(\Q)$ is dense in $\G(\A_\fin)$, where $\A_\fin$ is the ring of finite adeles.
Let $K_\fin$ be a given compact open subgroup of $\G(\A_\fin)$ and let $\Ga=\CG\cap K_\fin$.
Then $(\CG,\Ga)$ is a Hecke-pair and $\Ga$ can be chosen to be effective.
In this case the congruence completion is
$
\ol\CG\ \cong\ \G(\A_\fin).
$
Let $\widehat{\G(\A)}$ be the unitary dual of $\G(\A)$.
Note that every $\pi\in\widehat{\G(\A)}$ is a tensor product
$$
\pi\=\(\bigotimes_{p}\pi_p\)\ \otimes\ \pi_\infty,
$$
where the product runs over all primes $p$ and $\pi_p\in\widehat{\G(\Q_p)}$.
We also denote the representation $\bigotimes_{p}\pi_p$ of the group $\G(\A_\fin)$ by $\pi_\fin$.

Let $(\pi,V)$ be a representation of the locally compact group $\G(\A_\fin)$.
A vector $v\in V$ is called \e{smooth}, if $v$ is stabilized by an open subgroup of $\G(\A_\fin)$.
For any continuous representation, the set $V^\infty$ of smooth vectors is a dense, $\G(\A_\fin)$-stable subspace.

In this paper we assume that
$$
\Q-{\rm rank}(\G)\=0.
$$
this implies that $\G(\Q)\bs\G(\A)$ is compact.
For $\pi\in\widehat{\G(\A)}$ let
$$
m(\pi)\=\dim\Hom_{\G(\A)}(V_\pi,L^2(\G(\Q)\bs\G(\A))).
$$
Then $m(\pi)\in\N_0$ and 
$$
L^2(\G(\Q)\bs\G(\A))\ \cong\ \bigoplus_{\pi\in\widehat{\G(\A)}}m(\pi)\pi,
$$
Fix a congruence subgroup $\Sigma$ of $\Ga$.
Then there exists a compact open subgroup $K_\Sigma$ of $\G(\A_\fin)$ such that $\Sigma= K_\Sigma\cap \G(\Q)$.
The $G$-equivariant identification $\G(\Q)\bs\G(\A)/K_\Sigma\cong \Sigma\bs G$ gives an isomorphism of unitary $G$-representations
$$
L^2(\G(\Q)\bs\G(\A))^{K_\Sigma}\ \cong\ L^2(\Sigma\bs G).
$$
Hence we get
$$
L^2(\Sigma\bs G)\cong\bigoplus_{\pi_\infty\in\widehat G}m_{\Sigma}(\pi_\infty)\pi_\infty,
$$
where
$$
m_{\Sigma}(\pi_\infty)\=\sum_{\pi_\fin\in\widehat{\G(\A_\fin)}}m(\pi_\fin\otimes\pi_\infty)\,\dim(\pi_\fin^{K_\Sigma}).
$$

For $\pi\in \what G$ we write
$$
h^1(\pi)\=\dim\H^1(\g,K,\pi^\infty)\ \in\ \N_0.
$$
We have shown that for each congruence subgroup $\Sigma$ there is a natural inclusion
$$
\H_{q}^1(\Sigma,\C)\ \hookrightarrow\ \(\bigoplus_{\pi\in\what{\G(\A)}}m(\pi)h^1(\pi)\pi_\fin^{K_\Sigma}\)^{\otimes q}.
$$
Taking limits we see that there is a natural map of $\G(\A_\fin)$-representations
$$
\what\H_{q}^1(\C)\ \hookrightarrow\ \(\bigoplus_{\pi\in\what{\G(\A)}}m(\pi)h^1(\pi)\pi_\fin\)^{\otimes q}.
$$
We now state the main theorem.

\begin{theorem}
\begin{enumerate}[\rm (a)]
\item Let $q\ge 1$.
The $\G(\A_\fin)$-representation on $\what\H_{q}^1(\C)$ is the space of smooth vectors of a unitary representation $\what\H_{q}^1(\C)^u$ which is a subrepresentation of the tensor power
$$
\(\bigoplus_{\pi\in\what{\G(\A)}}m(\pi)h^1(\pi)\pi_\fin\)^{\otimes q}\ \cong\ \H^1(\g,K,L^2(\G(\Q)\bs\G(\A)))^{\otimes q}.
$$
\item
The $\G(\A)$-representation on $\what{\ol H}_q^0(L_\loc^2(G))$ is the space of $\G(\A_\fin)$-smooth vectors of a unitary representation which is isomorphic to
$$
\what\H_q^1(\C)^u\otimes L^2(\G(\Q)\bs \G(\A)).
$$
Here $\G(\A_\fin)$ acts on both tensor factors, but $G=\G(\R)$ acts only on the second.
\end{enumerate}
\end{theorem}

\prf Part (a) has been shown in this section above and part (b) is then a direct application of Proposition \ref{Prop3.4}.
\qed

{\small Mathematisches Institut\\
Auf der Morgenstelle 10\\
72076 T\"ubingen\\
Germany\\
\tt deitmar@uni-tuebingen.de}

\today

\end{document}